\documentclass[12pt,leqno]{article}

\usepackage{amssymb}
\usepackage{amsmath}
\usepackage{graphicx}
\usepackage[usenames]{color}

\def\PI{{\Pi}}

\def\hf{~\hfill\rule{1,5mm}{1,5mm}}

\def\derp#1#2{\displaystyle\frac{\partial#1}{\partial#2}}
\def\derpp#1{\displaystyle\frac{\partial}{\partial#1}}

\def\calf{{\mathcal{F}}}
\def\calo{{\mathcal{O}}}

\def\calq{{\mathcal{Q}}}

\def\E{{\mathbb{E}}}
\def\rr{{\mathbb{R}}}
\def\nn{{\mathbb{N}}}

\def\9{{\infty}}

\def\a{{\alpha}}
\def\b{{\beta}}

\def\g{{\gamma}}
\def\wt{\widetilde}
\def\ov{\overline}

\def\vf{{\varphi}}
\def\oo{{\omega}}
\def\ooo{{\Omega}}
\def\pp{{\partial}}
\def\D{{\Delta}}
\def\vp{{\varepsilon}}

\def\derp#1#2{\displaystyle\frac{\partial#1}{\partial#2}}
\def\barr{\begin{array}}
\def\earr{\end{array}}

\def\dd{\displaystyle}

\def\n{\noindent }
\def\pas{\mathbb{P}\mbox{-a.s.}}

\def\eq{equa\-tion}

\def\dd{\displaystyle}
\def\vsp{\vspace*{2mm}\\ }

\def\ff{\forall }

\def\({\left(}
\def\){\right)}
\def\<{\left<}
\def\>{\right>}

\newtheorem{theorem}{Theorem}[section]

\newtheorem{remark}[theorem]{Remark}

\title{Backward uniqueness of stochastic parabolic like equations driven by Gaussian multiplicative noise}
\author{Viorel Barbu\thanks{Al.I.
Cuza University and Octav Mayer Institute of Mathematics of the Romanian Academy} \and Michael R\"ockner\thanks{Fakult\"at f\"ur Mathematik, Universit\"at Bielefeld, D-33501 Bielefeld, Germany}}
\date{}

\begin{document}

\maketitle

\begin{abstract}   One proves here the backward uniqueness of solutions to stochastic semi\-linear parabolic equations and  also for the  tamed Navier--Stokes equations driven by linearly multiplicative Gaussian noises. Ap\-pli\-ca\-tions to approximate controllability of nonlinear stochastic parabolic equations with initial controllers are given. The method of   proof relies on the logarithmic convexity property known to hold for  solutions to li\-near evolution equations in Hilbert spaces with self-adjoint prin\-cipal~part.\medskip\\
{\bf Keywords:} stochastic parabolic equation, backward uniqueness, approximating controllability.\\
{\bf MSC:} 60H15; 47H05; 47J05.
\end{abstract}

\section{Introduction}

Consider the stochastic parabolic \eq
\begin{equation}\label{e1.1}
\barr{l}
dX(t)-\dd\sum^d_{i,j=1}\derpp{\xi_i}\,\(a_{ij}(t,\xi)\,\derp {X(t)}{\xi_j}\)dt+b(t,\xi)\cdot\nabla X(t)\,dt\vsp
\qquad\qquad+\psi(t,\xi,X(t))dt=X(t)dW(t)\mbox{ in }(0,T)\times\calo,\vsp
X(0,\xi)=x(\xi),\ \xi\in\calo;\qquad X(t,\xi)=0\mbox{ on }(0,T)\times\pp\calo,\earr\end{equation}where
$\calo\subset\rr^d,$ $1\le d<\9$, is a bounded and open domain with the smooth boundary $\pp\calo$ and $W$ is a Wiener process of the form
\begin{equation}\label{e1.2}
 W(t,\xi)=\sum^\9_{j=1}\mu_j e_j(\xi)\b_j(t),\ \xi\in\ov\calo,\ t\ge0.\end{equation}
Here $\{e_j\}^N_{j=1}\subset C^2(\ov\calo)$ is an orthonormal basis in $L^2(\calo)$, $\{\b_j\}^\9_{j=1}$ is an independent system of real-valued Brownian motions on a probability space $(\ooo,\calf,\mathbb{P})$ with the natural filtration$(\calf_t)_{t\ge0}$, and $\{\mu_j\}\subset\rr$ is such that
\begin{equation}\label{e1.2a}
\sum^\9_{j=1}\mu^2_j\|e_j\|^2_{C^2_b}<\9,\end{equation}
where $\|\cdot\|_{C^2_b}$ denotes the supnorm of the functions and its first and second order derivatives.

As regards the functions $a_{ij}:[0,T]\times\ov\calo\to\rr,$ $b:[0,T]\times\ov\calo\to\rr$ and $\psi:[0,T]\times\ov\calo\to\rr$, we assume that the following conditions hold
\begin{equation}\label{e1.3}
\barr{l}
a_{ij}\in C([0,T]\times\ov\calo),\derpp  t\ a_{ij}\in C([0,T]\times\ov\calo),\ \derpp{\xi_j}\,a_{ij}\in C([0,T]\times\ov\calo),\vsp
a_{ij}=a_{ji},\ \ff i,j=1,...,d,\ b\in C([0,T]\times\ov\calo;\rr^d),\ {\rm div}_\xi b\in C([0,T]\times\ov\calo),\earr\hspace*{-6mm}\end{equation}
\vspace*{-4mm}\begin{equation}\label{e1.4}
\barr{l}
\dd\sum^d_{i,j=1}a_{ij}(t,\xi)u_iu_j\ge\g|u|^2_d,\  \ff u=(u_1,...,u_d)\in\rr^d,\vsp (t,\xi)\in[0,T]\times\ov\calo,\earr\end{equation}
 where $\g>0$ and $|\cdot|_d$ is the Euclidean norm on $\rr^d$.
\begin{equation}\label{e1.5}
\barr{l}
\psi,\psi_\vp\in C([0,T]\times\ov\calo\times\rr),\ \psi(t,\xi,0)\equiv0,\vsp
|\psi_\vp(t,\xi,r)|\le C(1+|r|_d),\ \ff (t,\xi,r)\in[0,T]\times\ov\calo\times R.\earr\end{equation}Moreover, $r\to\psi(t,\xi,r)$ is monotonically nondecreasing and
\begin{equation}\label{e1.6}
\barr{r}
|\psi(t,\xi,r_1)-\psi(t,\xi,r_2)|\le L|r_1-r_2||\psi_0(t,\xi,r_1,r_2)|,\vsp \ff r_1,r_2\in\rr,\ (t,\xi)\in[0,T]\times\ov\calo,\earr\end{equation}where $\psi_0\in C([0,T]\times\ov\calo\times\rr\times\rr)$ and $L>0$,
\begin{equation}\label{e1.7}
|\psi_0(t,\xi,r_1,r_2)|\le C(|r_1|^q+|r_2|^q+1),\ \ff r_1,r_2\in\rr,\ (t,\xi)\in[0,T]\times\ov\calo,\end{equation}where
\begin{equation}\label{e1.8}
\barr{ll}
0\le q<\dd\frac{d+2}{d-2}&\mbox{ if }d>2,\vsp
q\in(1,\9)&\mbox{ if }d=2,\earr\end{equation}and no polynomial growth condition of the form \eqref{e1.8} is necessary if $d=1$.

In the following, we denote by $L^2(\calo)$ the space of Lebesgue square integrable functions on $\calo$ with the norm denoted  by  $|\cdot|_2$ and the scalar product $\<\cdot,\cdot\>$.
We denote by $W^{m,p}(\calo)$, $H^1_0(\calo)$ and $H^{-1}(\calo)$ the standard Sobolev spaces on $\calo$ with the usual norms $\|u\|_{m,p}$, $\|\cdot\|_1$ and $\|\cdot\|_{-1}$, respectively.

We note that, under assumptions \eqref{e1.2}--\eqref{e1.8}, for each $x\in L^2(\calo)$, equation  \eqref{e1.1} has a unique solution $X$ satisfying
\begin{eqnarray}
&X\in L^\9(0,T;L^2(\ooo;L^2(\calo)))\cap L^2([0,T)\times\ooo;H^1_0(\calo))\label{e1.9}\\[2mm]
&\E\dd\int^T_0\left\|e^{W(t)}\ \frac d{dt}\,(e^{-W(t)}X(t))\right\|^2_{-1}dt<\infty.\label{e1.10}
\end{eqnarray}Moreover, if $x\in L^\9(\calo)$, then $X\in L^\9((0,T)\times\calo),$ $\pas$ (See \cite{2}, Corollary 6.1 and \cite{3aa}, Theorem 2.1, p. 425.)
By a solution to \eqref{e1.1}, we mean an $(\calf_t)_{t\ge0}$--adapted process  $X:[0,T]\to L^2(\calo)$, with continuous sample paths, which satisfies the equation
\begin{equation}\label{e1.11}
\barr{l}
X(t,\xi)-\dd\int^t_0\sum^d_{i,j=1}\derpp{\xi_i}\,
\left(a_{ij}(s,\xi)\,\derp X{\xi_j}\,(s,\xi)\right)ds\\
\qquad+\dd\int^t_0(b(s,\xi)\cdot\nabla X(s,\xi)+\psi(s,\xi,X(s,\xi)))ds\\
\qquad=x(\xi)+\dd\int^t_0 X(s,\xi)dW(s),\ t\in[0,T],\ \xi\in\calo,\ \pas\earr\end{equation}
(Here, $\derpp{\xi_i}$ are taken in sense of distributions.)

For deterministic linear parabolic equations of the form \eqref{e1.1} (with $\psi\equiv0)$ and, more generally, for linear evolution equations in Hilbert spaces with principal part self-adjoint  of class $C^1$ with respect to $t$, it  is known that one has backward uniqueness of solutions $X$, that is, if $X_1(T)=X_2(T)$, then $X_1\equiv X_2$. (See \cite{1}, \cite{3}, \cite{4}.) Here, we shall prove that such a result remains true in the stochastic case \eqref{e1.1}. A few consequences of this result to approximate controllability with respect to the initial data $x$ are derived and, in particular, the approximate controllability of \eqref{e1.1} with respect to the initial data $x$. In Section 4, we  prove a similar result for   tamed stochastic Navier--Stokes equations.

Other results in this direction were obtained recently in \cite{5az}. However, only for linear SPDE,  the method used here is completely different. In particular, in contrast to \cite{5az}, we obtain pathwise estimates, instead of estimates in expectation.

The idea of the proof in the parabolic case is to reduce equation \eqref{e1.1} by  a rescaling  procedure   to a random parabolic equation and apply to this equation the standard calculation to prove that $\log|y_1-y_2|^2$ is quasi-concave in $t$, where $y_i=e^{-W}X_i$. This procedure allows us to obtain sharp estimates on \mbox{$|X_1(t)-X_2(t)|_2$} as a function of $|X_1(T)-X_2(T)|_2$. The rescaling approach can also be done for $2D$ stochastic Navier--Stokes equations (see \cite{3prim}), but  only for (analytically) weak solutions  which have not enough regularity to apply the above arguments. As a first step, we therefore consider stochastic tamed  Navier--Stokes equations (see \cite{11prim}, \cite{6}, \cite{7}) in this paper. But, in this case,  we use a direct approach based on a computation of \mbox{$d(\log|X_1(t)-X_2(t)|^2)$} via It\^o's formula, which still leads  to the backward uniqueness, but the obtained estimates are only in expectation. As a matter of fact, the advantage of a rescaling procedure is that it provides pathwise estimates. Its implementation is, however, much harder for stochastic tamed Navier--Stokes equations.

\section{The first main result}
\setcounter{equation}{0}

Everywhere in the following, we assume that conditions \eqref{e1.2a}--\eqref{e1.8} are sa\-tisfied.

\begin{theorem}\label{t1} Let $X_1,X_2$ be two solutions to \eqref{e1.1}, such that $X_1(0),X_2(0)\in L^\9(\calo)$, $\pas$ Then, if $X_1(T)\equiv X_2(T)$, $\pas$, we have $X_1\equiv X_2$. Moreover, there is a random variable $\g^*:\ooo\to\rr$,  such that $\pas$,
\begin{equation}\label{e2.1}
|X_1(t)-X_2(t)|_2\le\exp
\(\frac{\g^*\|X_1(t_0)-X_2(t_0)\|^2_1}
{|X_1(t_0)-X_2(t_0)|^2_2}\)|X_1(T)-X_2(T)|_2,\end{equation}for all $t\in[t_0,T],$  where $t_0$ is arbitrary in $(0,T).$

If the function $r\to\psi(t,\xi,r)$ is Lipschitz $($uniformly in $(t,\xi))$, then \eqref{e2.1} extends to all $X_1,X_2$ with $X_1(0),X_2(0)\in L^2(\calo)$.\end{theorem}

As will explicitly be seen in the proof, $\g^*$ is given by \eqref{e3.9a}, \eqref{e3.9aa} and it depends on $(\|e^{-W}X_1\|^q_{L^\9(\calq)}
+\|e^{-W}X_2\|^q_{L^\9(\calq)})$ and $W$.  Here, $\calq=(0,T)\times \calo$.

As a direct consequence of Theorem \ref{t1}, we obtain the following approxi\-mate controllability result.

\begin{theorem}\label{t2} Assume further that $r\to\psi(t,\xi,r)$ is differentiable and that $\derpp r\,\psi\in L^\9((0,T)\times\calo\times\rr)$. Then $\pas$ the space $\{X^x(T);\ x\in L^2(\calo)\}$ is dense in $L^2(\calo)$. $($Here, $X^x$ is the solution to \eqref{e1.1}.$)$\end{theorem}

 In the control theory literature, this property is called the {\it approximate controllability with respect to the initial data $x$ which is viewed as a start controller} (see, e.g., \cite{5}).

\section{Proofs}
\setcounter{equation}{0}

\subsection*{Proof of Theorem \ref{t1}}

By the transformation $X=e^Wy$, we reduce \eqref{e1.1} to the random parabolic equation
\begin{equation}\label{e3.1}
\barr{l}
\derp yt-e^{-W}\dd\sum^d_{i,j=1}\,\derpp{\xi_i}\,\(a_{ij}\,\derpp{\xi_j}\,(e^Wy)\)+\mu y+e^{-W}\psi(t,\xi,e^Wy)=0\vspace*{-2mm}\\\hfill\mbox{ in }(0,T)\times\calo,\\
y(0,\xi)=x(\xi),\ \xi\in\calo;\qquad y=0\mbox{ on }(0,T)\times\pp\calo,\earr\end{equation}where
$$ \mu(\xi)=\dd\frac12\,\sum^\infty_{j=1}\mu^2_je^2_j(\xi),\ \xi\in\ov\calo,\ t\in[0,T].$$Equivalently,
\begin{equation}\label{e3.2}
\barr{l}
\derp yt-\sum^d_{i,j=1}\derpp{\xi_i}\,\(a_{ij}\,\derp y{\xi_j}\)+a_0y+a_1\cdot\nabla y+e^{-W}\psi(t,\xi,e^Wy)=0\vspace*{-2mm}\\\hfill\mbox{ in }(0,T)\times\calo,\\
y(0,\xi)=x(\xi),\ \xi\in\calo;\ \ y=0\mbox{ on }(0,T)\times\pp\calo,\earr\end{equation}where $a_0:[0,T]\times\ov\calo\to\rr,\ a_1:[0,T]\times\ov\calo\to\rr^d$ are given by

\begin{eqnarray}
& \dd a_0=\mu(\xi)+\sum^d_{i,j=1}a_{ij}\(\derp{^2W}{\xi_i\pp\xi_j}+\derp W{\xi_i}\,\derp W{\xi_j}\)+\sum^d_{i,j=1}\derp{a_{ij}}{\xi_i}\,\derp W{\xi_j}\label{e3.3}\\[2mm]
&\dd a_1=\left\{2\sum^d_{i=1}
a_{ij}\ \derp{W}{\xi_j}\right\}^d_{j=1}.
\label{e3.4}
\end{eqnarray}We refer to \cite{2} for a rigorous justification of this rescaling argument and for the equivalence of \eqref{e1.1} and \eqref{e3.1} as well as the precise formulation of the latter.

We set $H=L^2(\calo)$, $V=H^1_0(\calo),\ V^*=H^{-1}(\calo)$ with the norms $|\cdot|_2,$ $\|\cdot\|_1$ respectively $\|\cdot\|_{-1}$ and define the operators $A(t):V\to V^*$, $B(t):V\to H$ and $B_1(t):V\to V^*$, $t\in[0,T]$, by
$$\barr{rcl}
\<A(t)y,z\>&=&\dd\sum^d_{i,j=1}\int_\calo a_{ij}(t,\xi)\,\derp y{\xi_i}\,(t,\xi)\,
\derp z{\xi_j}\,(t,\xi)d\xi,\ \ff y,z\in V,\vspace*{2mm}\\
(B(t)y)(\xi)&=&a_0(t,\xi)y(\xi)+a_1(t,\xi)\cdot\nabla y(\xi),\ \xi\in\calo,\ y\in V,\vspace*{3mm}\\
(B_1(t)y)(\xi)&=&e^{-W(t,\xi)}\psi(t,\xi, e^{W(t,\xi)}y(\xi)),\ \xi\in\calo, y\in V.\earr$$
Here, $\<\cdot,\cdot\>$ is the pairing between $V$ and $V^*$ which  coincides with the scalar product of $H$ on $H\times V$.

We note that there exist $\alpha_1,\alpha_2>0$ such that for $t\in[0,T]$
\begin{eqnarray}
\|A(t)z\|_{-1} \le \alpha_1\|z\|_1,\ \ff z\in V,\label{e3.4az}\\[2mm]
\<A(t)z,z\> \ge \alpha_2\|z\|^2_1,\ \ff z\in V.\label{e3.4aaz}\end{eqnarray}
We rewrite \eqref{e3.2} as
\begin{equation}\label{e3.5}
\barr{l}
\dd\frac{dy}{dt}\,(t)+A(t)y(t)+B(t)y(t)+B_1(t)y(t)=0,\ t\in(0,T),\vsp
y(0)=x.\earr\end{equation}For $x\in L^2(\calo)$ and $\mathbb{P}$-a.e. $\oo\in\ooo$, equation \eqref{e3.2} (equivalently \eqref{e3.5})  has a unique solution
$$y\in C([0,T];L^2(\calo))\cap L^2(0,T;H^1_0(\calo)),\ \dd\frac{dy}{dt}\in L^2(0,T;H^{-1}(\calo)).$$(See, e.g., \cite{6az}.) By the smoothing effect of solutions on initial data we have also that
\begin{equation}\label{e3.7a}
Ay,\frac{dy}{dt}\in L^2(\delta,T;L^2(\calo)),\ y\in C([\delta,T];H^1_0(\calo)),\end{equation}for all $\delta\in(0,T).$
This follows by the following arguments.

Consider the approximating equation
\begin{equation}\label{e3.7aa}
\barr{l}
\dd\frac{dy_\vp(t)}{dt}+A(t)y_\vp(t)+B(t)y_\vp(t)+ B^\vp_1(t)y_\vp(t)=0,\vsp
y_\vp(t)=x,\earr\end{equation}
where $B^\vp_1(t)y_\vp(t)=e^{-W(t)}\psi_\vp(t,\cdot,e^{W(t)}y_\vp(t))$ and $\psi_\vp(r)=\psi((1+\vp\psi)^{-1}r)$,\break $\ff r\in\rr$, is the Yosida approximation of $r\to\psi(\cdot,r)$.

Multiplying \eqref{e3.7aa} by $y_\vp(t)$ and integrating over $(0,t)\times\calo$, we get
\begin{equation}\label{e3.9az}
|y_\vp(t)|^2_2+\int^t_0\|y_\vp(s)\|^2_1
\le C(1+|x|^2_2),\ t\in[0,T],\end{equation}
where $C$ is independent of $\vp$.

Clearly, $y_\vp\to y$ in $C([0,T];L^2(\calo))\cap L^2(0,T;H^1_0(\calo))$ for $\vp\to0$. Moreover, since $B(t)y_\vp(t)+B^\vp_1(t)y_\vp(t)=f_\vp\in L^2(0,T;H)$, we have that
$$\sqrt{t}\ \frac{dy_\vp}{dt},\ \sqrt{t}\,A(t)y_\vp\in L^2(0,T;H),\ \ff t\in(0,y).$$

Now, as easily seen by the definition of $A(t)$ and by \eqref{e1.5}, we have
\begin{equation}\label{e3.10az}
\barr{l}
\<A(t)y_\vp(t),B(t)y_\vp(t)\>
+\<A(t)y_\vp(t),B^\vp_1y_\vp(t)\>\vsp
\qquad
\ge-C(\|y_\vp(t)\|^2_1+
|A(t)y_\vp(t)|_2
\|y_\vp(t)\|_1+1)\vsp
\qquad
\ge-\dd\frac12\,|A(t)y_\vp(t)|^2_2
-C(\|y_\vp(t)\|^2_1+1),\ \ff t\in[0,T].\earr\end{equation}
where $C$ is independent of $\vp$.

Then, multiplying \eqref{e3.7aa} by $tA(t) y_\vp(t)$ and integrating over $(0,t)\times\calo$, we get after some calculation involving \eqref{e3.9az}, \eqref{e3.10az} that
$$\barr{r}
t\<Ay_\vp(t),y_\vp(t)\>+\dd\int^t_0|A(s)y_\vp(s)|^2ds
 \le C\(\dd\int^t_0 \|y_\vp(s)\|^2_1ds+1\) \le C(|x|^2_2+1),\vsp \ff t\in[0,T].\earr$$This yields
$$\int^t_0 s\(|A(s)y_\vp(s))^2+\left|\frac{dy_\vp}{dt}\,(s)\right|^2\)ds
\le C(1+|x|^2),$$where $C$ is independent of $\vp$. Then, \eqref{e3.7a} follows.

Moreover, if $x\in L^\9(\calo)$, then $y\in L^\9((0,T)\times\calo)$ (see \cite{7a}, Theorem 2.1, p. 425). It follows also that the process $t\to y(t)$ is $(\calf_t)_{t\ge0}$--adapted.

Let $y_1,y_2$ be two solutions to \eqref{e1.1} with $y_i(0)\in L^\9(\calo)$ and $y_i=e^{-W}X_i,$ $i=1,2$. We set
$$g(t,\xi)=\left\{\barr{cl}
\dd\frac{B_1(t)y_1(t,\xi)-B_1(t)y_2(t,\xi)}{y_1(t,\xi)-y_2
(t,\xi)}&\mbox{ on }[y_1\ne y_2]\vsp
0&\mbox{ on }[y_1=y_2].\earr\right.$$
We have by \eqref{e1.6} \eqref{e1.7},  that
$$|g(t,\xi)|\le
C(\|y_1\|^q_{L^\9(\calq)}+\|y_2\|^q_{L^\9(\calq)}+1),\ \ff(t,\xi)\in\calq:=[0,T]\times\calo.$$
Hence $g\in L^\9(0,T)\times\calo) $ and $C$ is independent of $\oo$.

We set $z=y_1-y_2$ and get by \eqref{e3.5} that
\begin{equation}\label{e3.6}
\frac{dz}{dt}+A(t)z+B(t)z+g(t,\xi)z=0,\ \ t\in(0,T).\end{equation}We have that $z\in L^2(0,T;V)$, $g\in L^\9((0,T)\times\calo)$, $B(t)z\in L^2(0,T;H)$, and $\frac{dz}{dt},A(t)z\in L^2(\delta,T;H)$ for each $\delta\in(0,T)$. Moreover, by \eqref{e1.3}, we have that
 $\(\frac d{dt}\,A(t)\)z(t)\in C([0,T];V^*).$ It follows also that, by \eqref{e1.3}, we have
 $$\left|\(\frac d{dt}\, A(t)\)z(t)\right|_{-1}\le C\|z(t)\|_1,\ \ \ff t\in(0,T).$$
 (Here and everywhere in the following, we shall denote by the same symbol $C$ several positive constants independent of $\oo$.)

Arguing as in \cite{3}, \cite{4}, we get by \eqref{e3.6} that $\pas$
\begin{equation}\label{e3.7}
\barr{l}
\dd\frac d{dt}\,(\<A(t)z(t),z(t)\>|z(t)|_2^{-2})\vsp
\qquad=\(2\<A(t)z(t),\dd\frac{dz}{dt}\,(t)\>+\<\(\dd\frac d{dt}\,A(t)\)z(t),z(t)\>\)|z(t)|^{-2}_2 \vsp
\qquad
-2\<\dd\frac{dz}{dt}\,(t),z(t)\>
\<A(t)z(t),z(t)\>|z(t)|^{-4}_2,\ t\in[0,T].
\earr\end{equation}
Of course, \eqref{e3.7} holds on a maximal interval $[0,T^*]$, where $z(t)\ne0$. By uniqueness of the solution to the linear Cauchy problem associated with \eqref{e3.6}, $z(t)=0$ on $[T_0,T]$. Hence, replacing if necessary $T$ by $T^*$ we may assume without any loss of generality that $z(t)\ne0$ for $t\in[0,T]$.
We set $f(t)=B(t)z(t)+g(t)z$ and, by \eqref{e3.7}  we have, a.e., $t\in(0,T),$
\begin{equation}\label{e3.8}
\barr{l}
\!\!\!\!\!\!\!\!\dd\frac d{dt}\,(\<A(t)z(t),z(t)\>|z(t)|^{-2}_2)\vsp
=-2\<A(t)z(t),A(t)z(t)+f(t)\>|z(t)|^{-2}_2\vsp
+\<\(\dd\frac d{dt}\, A(t)\)z(t),z(t)\>|z(t)|^{-2}_2\vsp
+\,2\<A(t)z(t)+f(t),z(t)\>\<A(t)z(t),z(t)\>|z(t)|^{-4}_2\vsp
\le C\|z(t)\|^2_1|z(t)|^{-2}_2- [2|A(t)z(t)|^2_2+2\<A(t)z(t),f(t)\>]|z(t)|^{-2}_2\vsp
+\,2[\<A(t)z(t),z(t)\>^2+
\<f(t),z(t)\>\<A(t)z(t),z(t)\>]|z(t)|^{-4}_2\vsp
=C\|z(t)\|^2_1 |z(t)|^{-2}_2-
2\left[\left|A(t)z(t)+\dd\frac12\,f(t)\right|^2_2
-\dd\frac14\,|f(t)|^2_2\right]|z(t)|^{-2}_2\\
+2\,\left[\<A(t)z(t)+\dd\frac12\,f(t),z(t)\>^2
-\dd\frac14\<f(t),z(t)\>^2\right]|z(t)|^{-4}_2\vsp
\le C\|z(t)\|^2_1|z(t)|^{-2}_2
+ |f(t)|^2_2|z(t)|^{-2}_2\vsp
\le C \alpha^{-1}_2\<A(t)z(t),z(t)\>|z(t)|^{-2}_2
+|f(t)|^2_2|z(t)|^{-2}_2.
\earr\end{equation}On the other hand, by \eqref{e3.3}, \eqref{e3.4} we have
$$|f(t)|^2_2\le C_2(\nu_1+\g_2)^2\|z(t)\|^2_1,$$where
$$\barr{l}
\nu_1=\nu_1(\oo)
\dd=C\(\sup_{t\in[0,T]}\|W(t)\|_{C^2_b}
+\sup_{(t,\xi)\in\mathcal{Q}}
|\nabla W(t,\xi)|^2\)\vsp
 \g_2=\g_2(\oo)=
 \|y_1\|^q_{L^\9(\calq)}+\|y_2\|^q_{L^\9(\calq)}+1\vsp
 \quad\,
=\|e^{-W}X_1\|^q_{L^\9(\calq)}
 +\|e^{-W}X_2\|^q_{L^\9(\calq)}+1.\earr$$
 Then, substituting into \eqref{e3.8} yields, for $t\in(0,T),$
$$\dd\frac d{dt}\,(\<A(t)z(t),z(t)\>|z(t)|^{-2}_2)
\le(C_1+C_2(\nu_1+\g_2)^2
|\<A(t)z(t),z(t)\>\!|\,|z(t)|^{-2}_2,$$where $C_1,C_2$ are independent of $\oo$. Hence
\begin{equation}\label{e3.9}
\<A(t)z(t),z(t)\>|z(t)|^{-2}_2
\le\exp(\g^*_1(t-t_0))\<A(t_0)z(t_0),z(t_0)\>|z(t_0)|^{-2}_2,
\end{equation}
for $t_0\le t\le T$. Here, $\g^*_1$ is the random variable
\begin{equation}\label{e3.9a}
\g^*_1(\oo)=C_1+C_2
 (\nu_1(\oo)+\g_2(\oo))^2.\end{equation}

On the other hand, we have
$$\barr{l}
\dd\frac12\, \frac d{dt}\log(|z(t)|^2_2)=-\<A(t)z(t)+B(t)z(t)+g(t)z(t),z(t)\>|z(t)|^{-2}_2\vsp
\qquad\qquad\ge-\<A(t)z(t),z(t)\>|z(t)|^{-2}_2
-C_2(\nu_1+\g_2)\|z(t)\|_1|z(t)|^{-1}_2\vsp
\qquad\qquad\ge-C_3(\nu_1+\g_2+1)
\<A(t)z(t),z(t)\>|z(t)|_2^{-2},\ \ff t\in(0,T).\earr$$Then, by \eqref{e3.9}  we obtain
$$\barr{r}
\dd\frac12\,\frac d{dt}(\log|z(t)|^2_2)
\ge-C_3(\nu_1{+}\g_2{+}1)\exp(\g^*_1(t{-}t_0))
\<A(t_0)z(t_0),z(t_0)\>
|z(t_0)|^{-2}_2,\earr$$for all $0<t_0<t<T$.

Integrating from $t$ to $T$, we obtain estimate \eqref{e2.1}, where
\begin{equation}\label{e3.9aa}
\g^*=C_4(\nu_1+\g_2+1)\frac1{\g^*_1}
\exp(\g^*_1(T-t_0)),\end{equation}
where $C_4$ is independent of $\oo$.
If $\psi$ is Lipschitz in $r$ uniformly with respect to $(t,\xi)$, then $g\in L^\9(0,T)\times\calo)$ for all $X_i$ with $X_i(0)\in L^2(\calo)$, $i=1,2$, and so condition $X_i(0)\in L^\9(\calo)$ is no longer necessary.  This completes the proof.\hf

\subsection*{Proof of Theorem \ref{t2}}

Denote by $S(t):L^2(\calo)\to L^2(\calo)$ the flow
\begin{equation}\label{e3.10}
S(t)x=y^x(t),\ t\in(0,T),\end{equation}
where $y^x$ is the solution to \eqref{e3.5} (equivalently \eqref{e3.2}). It is  easily seen that $S(T)$ is Fr\'echet differentiable on $L^2(\calo)$ and its Fr\'echet derivative at the origin $\Gamma:L^2(\calo)\to L^2(\calo)$ is given by $\Gamma u=DS(T)(0)u=v(T)$, where $v$ is the solution to the equation
\begin{equation}\label{e3.11}
\barr{l}
\dd\frac{dv}{dt}+A(t)v+B(t)v+e^{-W}\psi_r(t,\xi,e^W\wt y)v=0\mbox{ in }(0,T),\vsp
v(0)=u,\earr\end{equation}where $\wt y$ is the solution to \eqref{e3.5} with $\wt y(0)=0$ and $\psi_r=\derpp r\,\psi.$

Then, the dual operator $\Gamma^*:L^2(\calo)\to L^2(\calo)$ is given by $\Gamma^*p=z(0)$, $\ff p\in L^2(\calo)$, where $z$ is the solution to backward dual equation
\begin{equation}\label{e3.12}
\barr{l}
\dd\frac{dz}{dt}-A(t)z-B^*(t)z-e^{-W(t)}\psi_r(t,\xi,e^{W(t)}\wt y(t))z=0,\ t\in(0,T),\vsp
z(T)=p,\earr\end{equation}which, clearly, is well posed for all $p\in L^2(\calo)$.

By Theorem \ref{t1} (applied to the backward equation \eqref{e3.12}), the operator $\Gamma^*$ is injective on $L^2(\calo)$ and, as well known (see e.g., Proposition IV.1 in \cite{3}), this implies that the space $\{S(T)x;\ x\in L^2(\calo)\}$ is dense in $L^2(\calo)$, as claimed. \hf

\begin{remark}\label{r3.1}\rm One might ask whether, under the assumptions of  Theorem \ref{t2}, the set  $\{X(T,x);\ x\in L^2(\calo)\}$ is dense in $L^2(\ooo;\,L^2(\calo))$, that is, in the mean square norm $(\E|\cdot|^2_2)^{\frac12}$. Clearly, this happens if the stochastic backward dual equation associated with \eqref{e1.1}, that is,
\begin{equation}\label{e3.13}
\left\{\barr{l}
dp+\dd\sum^d_{i,j=1}\derpp{\xi_i}\,\(a_{ij}(t,\xi)\,\derp p{\xi_j}\)-{\rm div}(bp)dt-\psi_r(t,\xi,0)pdt\vsp
\qquad+\dd\sum^\9_{j=1}\mu_je_jq_jdt
=\dd\sum^\9_{j=1}q_j(t)d\b_j(t)\mbox{ in }(0,T)\times\calo,\vsp
p=0\mbox{ on }(0,T)\times\pp\calo,\earr\right.\end{equation}
has the forward uniqueness property, that is, $p(0)\equiv0$ implies $p\equiv0
$. However, as far as we know, this is an open problem.
\end{remark}

\section{The second main result, the backward\\ uniqueness for stochastic  $3D$-tamed\\ Navier--Stokes equations}
\setcounter{equation}{0}

Consider the stochastic equation
\begin{equation}\label{e4.1}
\barr{l}
dX-\nu\Delta Xdt+(X\cdot\nabla)Xdt+g_N(|X|^2_3)Xdt=XdW+\nabla pdt\\\hfill\mbox{ in }(0,T)\times\calo,\vsp
\nabla\cdot X=0\mbox{ in }(0,T)\times\calo;\quad X=0\mbox{ on }(0,T)\times\pp\calo,\vsp
X(0)=x\mbox{ in }\calo,\earr\end{equation}where $\calo$ is a bounded and open subset of $\rr^3$,   with smooth boundary $\pp\calo$ and $|\cdot|_3$ denotes the Euclidean norm on $\rr^3$. $W$ is the Wiener process on $(\ooo,\calf,(\calf_t),\mathbb{P})$ from the previous sections, i.e.,
$$W(t,\xi)=\sum^\9_{j=1}\mu_je_j(\xi)\b_j(t),\ \ \xi\in\ov\calo,\ t\ge0,$$where $\{e_j\}\subset C^2(\ov\calo)$ is an orthonormal basis in $L^2(\calo)$, but with $\mu_j\in\rr$ satisfying the stronger condition
\begin{equation}\label{e4.2}
\sum^\9_{j=1}\mu^2_j(|e_j|^2_\9+|\nabla e_j|^2_\9)<\9,\end{equation}where $|\cdot|_\9$ is the norm in $L^\9(\calo)$.
Here, $g_N\in C^1(\rr^+)$, $N\in\nn$, is a given function such that
\begin{equation}\label{e4.2a}
g_N(r)=\left\{\barr{ll}
0&\mbox{ for }r\in[0,N],\vsp
\dd\frac{r-N-1}{\nu}&\mbox{ for }r\ge N+1,\vsp
0\le g'_N(r)\le C&\mbox{ for }r\in\rr.\earr\right.\end{equation}
Equation \eqref{e4.1} is a modified version of the stochastic Navier--Stokes systems and was introduced by R\"ockner and Zhang \cite{11prim} (see also \cite{7}, \cite{6}). In the deterministic case, any bounded solution to the standard stochastic Navier--Stokes equations is a solution to \eqref{e4.1} for sufficiently large $N$. However, in contrast to the case of the standard stochastic $3D$--Navier--Stokes equation, which in general has a (probabilistically) weak solution only (see, e.g., \cite{3a}, \cite{3aa}, \cite{5a}), problem \eqref{e4.1} is well posed in the (probabilistically) strong sense in an appropriate space, even in $3{-}D$.

\begin{remark}\rm In all what follows we could have taken a more general noise term than $XdW$, more precisely, the same type of noise as in \cite{11prim}. All the arguments are exactly the same in this more general case. However, we restrict ourselves to $XdW$ for simplicity and in order not to change the frame in comparison to Sections 2 and 3.\end{remark}

By strong solution to \eqref{e4.1}, we mean a pair of  $(\calf_t)_{t\ge0}$--adapted processes $X:[0,T]\to H=\{y\in(L^2(\calo))^3;\,\nabla\cdot y=0,\ y\cdot\vec n=0\mbox{ on }\pp\calo\}$, $p:[0,T]\to H^1(\calo)$, which has continuous sample paths and satisfies
$$\barr{l}
X\in L^\9(0,T;L^2(\ooo;H))\cap L^2((0,T)\times\ooo;(H^1_0(\calo))^3)\vsp
X(t)=\nu\dd\int^t_0\Delta X(s)ds-\int^t_0((X(s)\cdot\nabla)X(s)+g_N(|X(s)|^2_3)X(s))ds\vsp
\hfill+\dd\int^t_0\nabla p(s)ds+\int^t_0X(s)dW(s),\ \ff t\in[0,T],\ \pas\earr$$
For each $x\in (H^1_0(\calo))^3\cap H$, equation \eqref{e4.1} has a unique strong solution $X$, which satisfies
\begin{equation}\label{e4.3}
\E\left[\dd\sup_{0\le t\le T}\|X(t)\|^2_{(H^1_0(\calo))^3}\right]
+\dd\int^T_0\E\left[\|X(t)\|^2_{(H^2(\calo))^3}\right]dt\le C\|x\|^2_{(H^1_0(\calo))^3}.\end{equation}(See \cite{6}, Theorem 3.1.)

Furthermore, since our initial condition is not random, we also have (see \cite[Formula (3.12)]{6})
\begin{equation}\label{e4.4prim}
\int^T_0\E[\|X(t)\|^6_{(H^1(\calo))^3}]dt<\9.\end{equation}

In the following, we shall use the standard notations
$$\barr{l}
V=(H^1_0(\calo))^3\cap H,\ A=- \PI\
\Delta,\ D(A)=(H^2(\calo))^3\cap V,\earr$$
where $\PI:(L^2(\calo))^3\to H$ is the Leray projection (see  \cite{8}). We set, also,
$$b(y,z,\theta)=\int_\calo y_iD_iz_j\theta_jd\xi,\ \ \ff y,z,\theta\in V,$$and denote by $B:V\to V^*$ the operator $$\<BX,\vf\>=b(X,X,\vf),\ \ff \vf\in V.$$
The norm of $V$ will be taken as
$$\|y\|=\<Ay,y\>^{\frac12},$$
where $\<\cdot,\cdot\>$ is the duality pairing between $V$ and its dual $V^*$. This norm is equivalent to $\|\cdot\|_{(H^1(\calo))^3}.$ On $V\times H$, this is just the scalar product of $H$. The norm of $H$ is denoted in the following by $|\cdot|$. We recall that we have
\begin{equation}\label{e4.3a}
|b(y,z,\theta)|\le C\|y\|_{m_1}\|z\|_{m_2+1}\|\theta\|_{m_3},\end{equation}where $m_1+m_2+m_3\ge\frac 32$, if $m_i\ne\frac 32$, and $m_1+m_2+m_3>\frac 32$, if one of the $m_i$ is~$\frac 32$.
 (Here, $\|\cdot\|_m$ is the norm of the Sobolev space $H^m(\calo)$.)

Then, we can rewrite \eqref{e4.1} as the stochastic differential equation on the space $H$\vspace*{-3mm}
\begin{equation}\label{e4.4}
\barr{l}
dX+\nu AXdt+BXdt+\PI(g_N(|X|^2_3)X)dt=\dd\sum^\9_{j=1}\mu_j
\PI(Xe_j)d\b_j,\vspace*{-2mm}\\\hfill t\in(0,T),\vspace*{-2mm}\\
X(0)=x.\earr\end{equation}
It is useful for the time being to mention that, as shown in \mbox{\cite[Theorem 3.1]{6},}  the solution $X$ to \eqref{e4.4} is obtained as
\begin{equation}\label{4.6az}
\barr{ll}
X=\dd\lim_{n\to\9}u_n&\mbox{ weakly in }L^2(\ooo_T;(H^2(\calo))^3)\\
&\mbox{ weakly-star in }L^2(\ooo;L^\9(0,T;(H^1_0(\calo))^3)),\vsp
\Pi_nF(u_n)\to F&\mbox{ weakly in }L^2(\ooo_T;H),\earr\end{equation}
where $$\barr{lcl}
\ooo_T&=&[0,T]\times\ooo,\vsp
Fu&=&-\nu Au-Bu-\Pi(g_N(|u|^2_3)u)\earr$$and $\Pi_n$ is the orthogonal projection of $H$ onto $H_n={\rm span}(\wt e_1,\wt e_2,..,\wt e_n),$ $\{\wt e_i,\ i\ge1\}\subset (H^2(\calo))^3\cap V$ being a fixed orthonormal basis in $H$ consisting of eigenvectors of $A$. Moreover, $u_n$ is the solution to the finite dimensional stochastic differential equation
\begin{equation}\label{e4.6aaz}
\barr{lcl}
du_n(t)&=&\Pi_n F(u_n(t))dt
+\dd\sum^\9_{j=1}\mu_j\Pi_n(u_ne_j)d\beta_j,\\
u_n(0)&=&\Pi_n x.\earr\end{equation}

Theorem \ref{t3} below is the backward uniqueness result for strong solutions to \eqref{e4.1}.

\begin{theorem}\label{t3} Let $X_1,X_2$ be two solutions to \eqref{e4.1}, which satisfy \eqref{e4.3}. Then,   for any pair of solutions $X_1,X_2$ of \eqref{e4.1},
\begin{equation}\label{e4.4a}
\barr{r}
\E\left[e^{-C\g(t)}\log(|X_1(t){-}
X_2(t)|^2)\right]\le
\E\left[ e^{-C\g(T)}\log (|X_1(T){-}X_2(T)|^2)\right]\vsp
+\,C+\|X_1(0)-X_2(0)\|\,
|X_1(0)-X_2(0)|^{-1},\
 \ff t\in(0,T),\earr\hspace*{-3mm}\end{equation}
where
\begin{equation}\label{e4.6a}
\barr{r}
\g(t)=\dd\int^t_0
(\|X_1(s)\|^2_{{(W^{1,4}(\calo))^3}}+
\|X_2(s)\|^2_{{(W^{1,4}(\calo))^3}}\vsp
+\|X_1(s)\|^4+\|X(s)\|^4+1)ds,\  t\ge0,\earr\end{equation}which is finite $\pas$  by \eqref{e4.3}, \eqref{e4.4prim}, and $C$ is a positive constant independent of $\oo\in\ooo$. Furthermore, the last summand in \eqref{e4.4a} is defined to be zero, if $X_1(0)=X_2(0)$. In particular, in the deterministic case, i.e.  when the noise is zero, it
follows that $X_1(T) = X_2(T)$ implies $X_1(t) = X_2(t)$ for all $t\in [0,T]$.\end{theorem}

\begin{remark}\rm The expectations in \eqref{e4.4a} are well defined because of \eqref{e4.3}, but maybe equal to $-\9$, as happens in the case when $X_1(T)=X_2(T)$ $\pas$\end{remark}

 \n{\bf Proof of Theorem \ref{t3}.} For simplicity, we shall take  $\nu=1$ in the following.

  We set $Z=X_1-X_2$ and, by \eqref{e4.4}, we  get for $Z$ the linear equation
\begin{equation}\label{e4.3b}
dZ+AZdt+F_1dt+F_2dt=\PI(ZdW)\mbox{ in }(0,T),
\end{equation}
where $F_i:[0,T]\to H,\ i=1,2,$ are given by
\begin{eqnarray}
F_1&=&\PI((Z\cdot\nabla)X_1+(X_2\cdot\nabla)Z)\label{e4.4b}\\
F_2&=&\PI(g_N(|X_1|^2_3)X_1-g_N(|X_2|^2_3)X_2).\label{e4.5}
\end{eqnarray}
We have for $dt$-a.e. $t\in[0,T]$
\begin{equation}\label{e4.6}
\barr{ll}
|F_1(t)|\!\!\!&
\le C_1\|Z(t)\|(\|\nabla X_1(t)\|_{(L^{4}(\calo))^3}
+\|X_2(t)\|_{(L^\9(\calo))^3})\vsp
&\le C_2\|Z(t)\|(\|X_1(t)\|_{(W^{1,4}(\calo))^3}
+\|X_2(t)\|_{{(W^{1,4}}(\calo))^3}),\vsp
&\hfill \ff t\in[0,T],\earr
\end{equation}
because, by the Rellich--Kondrachev theorem (see, e.g., \cite[p.~285]{5a}),\break $(W^{1,4}(\calo))^3\subset L^\9(\calo)$.
By \eqref{e4.2a} it follows that
\begin{equation*}
\barr{r}
|g_N(|X_1|^2_3)X_1-g_N(|X_2|^2_3)X_2|_3
\le C_3(|X_1(t)|^2_3+|X_2(t)|^2_3)|Z(t)|_3\label{e4.8}\vsp \mbox{ a.e. in }(0,T)\times\calo\times\ooo,\earr
\end{equation*}which, by the Sobolev embedding, implies
\begin{equation}\label{e4.16x}
|F_2(t)|\le C_4(\|X_1(t)\|^2+\|X_2(t)\|^2)\|Z\|.\end{equation}
(Here and everywhere in the following,  $C_i$, $i=1,...,$ are positive constants independent of $\oo\in\ooo$.)

Also, in this case, we have (see \eqref{e4.6aaz})
\begin{equation}\label{e4.14az}
\barr{ll}
Z=\dd\lim_{n\to\9}z_n&\mbox{ weakly in }L^2(\ooo_T;(H^2(\calo))^3),\\
&\mbox{ weakly-star in }L^2(\ooo;L^\9(0,T;(H^1_0(\calo))^3)),\earr\end{equation}where
\begin{equation}\label{e4.14aaz} \barr{c}
dz_n+A_nz_ndt
+F^n_1dt+F^n_2dt
=\dd\sum^\9_{j=1}\mu_j\Pi_n(z_n e_j)d\beta_j,\vsp
z_n(0)=\Pi_n(Z(0)),\earr\end{equation}
where $F^n_i=\Pi_nF_i,\ i=1,2,$ $A_n=\Pi_nA.$
Moreover, estimates \eqref{e4.6}--\eqref{e4.8} hold in this case for $z_n$, $u^1_n=\Pi_n X_1,$ $u^2_n=\Pi_nX_2$ instead of $Z$, $X_1$ and $X_2$, respectively.

We consider the function
  $$\vf_\vp(y)=\frac{\|y\|^2}{|y|^2+\vp},\ \ y\in V,$$  where $\vp>0$ is arbitrary but fixed. We see that  $\vf_\vp$ is $C^2$ on $V$ and its Gateaux derivative $D\vf_\vp\in V'$ is given by
\begin{equation}\label{e4.9}
D\vf_\vp(y)=2[Ay(|y|^2+\vp)-y\|y\|^2](|y|^2+\vp)^{-2},\ y\in V.\end{equation}
Moreover, we have, for the second derivative  $D^2$,
\begin{equation}\label{e4.10}
\barr{ll}
D^2\vf_\vp(y)(h)
 {=}2(|y|^2{+}\vp)^{-2}
[(|y|^2{+}\vp)Ah{-}h\|y\|^2 {+}2Ay\<y,h\>{-}2y\<Ah,y\>]
\vsp
 -\,4\<y,h\>[Ay(|y|^2+\vp)-y\|y\|^2](|y|^2+\vp)^{-3},\ \ff y,h\in V.\earr\hspace*{-5mm}\end{equation}
If we heuristically apply It\^o's formula to $\vf_\vp$ in equation \eqref{e4.3b}, we get
\begin{equation}\label{e4.11}
\barr{l}
d\vf_\vp(Z(t))+2|AZ(t)|^2(|Z(t)|^2+\vp)^{-1}dt\vsp
\qquad
-\,2\|Z(t)\|^4(|Z(t)|^2+\vp)^{-2}dt\vsp
\qquad +\,2\<AZ(t),F_1(t)+
F_2(t)\>(|Z(t)|^2+\vp)^{-1}dt\vsp
\qquad-\,2\<Z(t),F_1(t)+F_2(t)\>\|Z(t)\|^2(|Z(t)|^2+\vp)^{-2}dt\vsp
\qquad =\,\dd\frac12\sum^\9_{j=1}\mu^2_j
\<D^2\vf_\vp(Z(t))(Z(t)e_j),Z(t)e_j\>dt \vsp
\qquad+\,\<D\vf_\vp(Z(t)),Z(t)dW(t)\>,\earr\hspace*{-10mm}
\end{equation}for $t\in[0,T]$.
However, it should be said that, since $Z$ is not a semi\-martingale in $V$, the It\^o formula cannot be applied in \eqref{e4.3b} and so to get \eqref{e4.11} we shall invoke a more sophisticated argument based on the appro\-xi\-ma\-ting equation \eqref{e4.14az}. Namely, we shall apply It\^o's formula in \eqref{e4.14aaz} to the function
$$\vf_\vp(v)=\|v\|^2\rho_\vp(|v|^2),\ \ff v\in V,$$
where $\rho_\vp(r)=\frac1{r+\vp},$ $\ff r\ge0$. For $\vf_\vp$, \eqref{e4.9} and \eqref{e4.10} remain valid, and so we~get
\begin{equation}\label{e4.17az}
\hspace*{-9mm}\barr{l}
d\vf_\vp(z_n(t))+2|A_nz_n(t)|^2
(|z_n(t)|^2+\vp)^{-1}dt\vsp
\qquad-2\|z_n(t)\|^4
(|z_n(t)|^2+\vp)^{-2}dt\vsp
\qquad+2\<A_nz_n(t),F^n_1(t)+F^n_2(t)\>
(|z_n(t)|^2+\vp)^{-1}dt\vsp
\qquad-2\<z_n(t),F^n_1(t)+F^n_2(t)\>
\|z_n(t)\|^2(|z_n(t)|^2+\vp)^{-2}dt\vsp
\qquad=\dd\frac12\sum^\9_{j=1}\mu^2_j
\<D^2\vf_\vp(z_n(t))\>
\<z_n(t)e_j,z_n(t)e_j\>dt\\
\qquad+\dd\sum^\9_{j=1}\<D \vf_\vp(z_n(t)),
z_n(t)dW(t)\>.\earr\hspace*{-10mm}\end{equation}
Taking into account \eqref{e4.14az}, we may pass to the limit in \eqref{e4.17az} and get  for $Z$ just formula \eqref{e4.11}.

We have
\begin{equation}\label{e4.12}
\barr{l}
|AZ|^2(|Z|^2+\vp)^{-1}-\|Z\|^4(|Z|^2+\vp)^{-2}\vsp
\qquad
-\,\<Z,F_1+F_2\>\|Z\|^2(|Z|^2+\vp)^{-2}
+\<AZ,F_1+F_2\>(|Z|^2+\vp)^{-1}\vsp
\qquad=(|Z|^2+\vp)^{-2}[|AZ|^2(|Z|^2+\vp)-\<AZ,Z\>^2\vsp
\qquad
+\,\<AZ,F_1+F_2\>(|Z|^2+\vp)-\<Z,F_1+F_2\>\|Z\|^2]\vsp
\qquad=(|Z|^2+\vp)^{-2}
\left[\left|AZ+\dd\frac12\,(F_1+F_2)\right|^2|Z|^2
-\dd\frac14\,|F_1+F_2|^2|Z|^2\right.\vsp
\qquad\left.-\(\<AZ,Z\>+\<\dd\frac12\,(F_1+F_2),Z\>\)^2
+\dd\frac14\<F_1+F_2,Z\>^2\right]\vsp
\qquad+\vp(|Z|^2+\vp)^{-2}(|AZ|^2+\<AZ,F_1+F_2\>)\vsp
\qquad
\ge-\dd\frac14\,(|Z|^2+\vp)^{-2}|Z|^2|F_1+F_2|^2\vsp
\qquad
+\ \vp(|Z|^2+\vp)^{-2}(|AZ|^2+\<AZ,F_1+F_2\>)\vsp
\qquad
\ge-\dd\frac14\,(|Z|^2+\vp)^{-1}|F_1+F_2|^2
-\frac\vp4\,|F_1+F_2|^2(|Z|^2+\vp)^{-2}\vsp
\qquad\ge-\dd\frac12\,(|Z|^2+\vp)^{-1}|F_1+F_2|^2.\earr\end{equation}
We have also by \eqref{e4.2} and \eqref{e4.10}
\begin{equation}\label{e4.13}
\barr{lcl}
\left|\dd\sum^\9_{j=1}\mu^2_j
\<D^2\vf_\vp(Z)(Ze_j),Ze_j\>\right|\vsp
\qquad\le C_4(|Z|^2+\vp)^{-1}\|Z\|^2
\dd\sum^\9_{j=1}\mu^2_j(|e_j|^2_\9+|\nabla e_j|^2_\9)
\le C_5\vf_\vp(Z).\earr\end{equation}
On the other hand, by \eqref{e4.6}, \eqref{e4.8} we see that
$$\!\!\barr{l}
|F_1+F_2|^2(|Z|^2+\vp)^{-1}\vsp
 \le C_6(\|X_1\|^2_{(W^{1,4}(\calo))^3}+
\|X_2\|^2_{(W^{1,4}(\calo))^3}
+\|X_1\|^4+\|X_2\|^4)
\|Z\|^2(|Z|^2+\vp)^{-1}\vsp
\le C_6(\g'-1)\vf_\vp(Z), \earr$$where $\g'$ is the derivative of $\g$ given in \eqref{e4.6a}.

Substituting \eqref{e4.12}, \eqref{e4.13} into \eqref{e4.11}, we obtain    that
\begin{equation}\label{e4.18}
d\vf_\vp(Z(t))
\le C_7\g'(t)\vf_\vp(Z(t))dt
 + \<D\vf_\vp(Z(t)),Z(t)dW(t)\>.\end{equation}
 We note that \eqref{e4.4prim} ensures the integrability of the integrands in the right hand side.

Integrating \eqref{e4.18} from $0$ to $t$, multiplying by $\exp(-C_7\g(t))$ and applying It\^o's product rule, we obtain
$$\barr{l}
e^{-C_7\g(t)}\vf_\vp(Z(t))\le
\vf_\vp(Z(0))
+\dd\int^t_0 e^{-C_7\g(s)}
\<D\vf_\vp(Z(s)),Z(s)dW(s)\>,\earr$$
and this yields
\begin{equation}\label{e4.19}
\E[\vf_\vp(Z(t))\exp(-C_7\g(t))]\le
\vf_\vp(Z(0)),\ \ff t\in[0,T].\end{equation}

Next, we apply the It\^o formula to \eqref{e4.3b} and the function $$\psi_\vp(z)=\frac12\log(|z|^2+\vp),\ z\in V.$$
Taking into account that $D\psi_\vp(z)=z(|z|^2+\vp)^{-1},$ we obtain that
$$\barr{r}
d\(\dd\frac12\,\log(|Z(t)|^2+\vp)\)+\vf_\vp(Z(t))dt
+\<F_1(t)+F_2(t),Z(t)\>(|Z(t)|^2+\vp)^{-1}dt\vsp
=\dd\frac12\sum^\9_{j=1}\mu^2_j
\<D^2\psi_\vp(Z(t)) (Z(t)e_j),Z(t)e_j\> dt+\<D\psi_\vp(Z(t)),Z(t)dW(t)\>.\earr$$
By \eqref{e4.6}, \eqref{e4.16x},
  we get after some calculations that
$$\barr{l}
d\(\dd\frac12\,\log(|Z(t)|^2+\vp)\)+\vf_\vp(Z(t))dt\vsp
\ge-C_8
(\|X_1(t)\|_{(W^{1,4}(\calo))^3}
+\|X_2(t)\|_{(W^{1,4}(\calo))^3}\vsp
+\|X_1(t)\|^2+\|X_2(t)\|^2)
|Z(t)|\|Z(t)\|
(|Z(t)|^2+\vp)^{-1}dt\vsp
-C_9dt  +\<D\psi_\vp(Z(t)),Z(t)dW(t)\>\vsp
\ge-\vf_\vp(Z(t))dt-C_{10}\g'(t)dt
+\<D\psi_\vp(Z(t)),Z(t)dW(t)\>,\earr$$
because, by \eqref{e4.2},
$$\sum^\9_{j=1}\mu^2_j\<D^2\psi_\vp(Z(t))(Z(t)e_j),Z(t)e_j\>\le C_9.$$
This yields
$$d(\log(|Z(t)|^2+\vp))
\ge
-4\vf_\vp(Z(t))dt-
C_{11}\g'(t)dt
+\<D\psi_\vp(Z(t)),Z(t)dW\>.$$
Letting $T\ge r>t$ and integrating, we obtain
$$\barr{ll}
\log(|Z(r)|^2+\vp) \ge\!\!\!&
\log(|Z(t)|^2+\vp)
-4\dd\int^r_t\vf_\vp(Z(s))ds\vsp
&-\,C_{11}\dd\int^r_t\g'(s)ds+\dd\int^r_t\<D\psi_\vp(Z(s)),Z(s)dW(s)\>.\earr$$Then, multiplying by $\exp(-C_{11}\g(t))$ and using It\^o's product formula, we get as above
$$\barr{l}
e^{-C_{11}\g(r)}\log(|Z(r)|^2+\vp)
=e^{-C_{11}\g(t)}\log(|Z(t)|^2+\vp)\vsp
+\dd\int^r_t e^{-C_{11}\g(s)}d\log(|Z(s)|^2+\vp)
-C_{11}\dd\int^r_t\log(|Z(s)|^2+\vp)\g'(s)e^{-C_{11}\g(s)}ds\vsp
\ge e^{-C_{11}\g(t)}\log(|Z(t)|^2+\vp)-\dd\int^r_t e^{-C_{11}\g(s)}
(4\vf_\vp(Z(s))+C_{11}\g'(s))ds\vsp
+\dd\int^r_t e^{-C_{11}\g(s)}
\<D\psi_\vp(Z(s)),X(s)dW(s)\>\vsp
-C_{11}\dd\int^r_t\log(|Z(s)|^2+\vp)\g'(s)e^{-C_{11}\g(s)}ds.\earr$$
Taking $r=T$, we get taking expectation
$$\barr{l}
\E[e^{-C_{11}\g(t)}\log(|Z(t)|^2+\vp)]
\le \E[e^{-C_{11}\g(T)}\log(|Z(T)|^2+\vp)]\vsp
\qquad+\,4
\E\dd\int^T_t e^{-C_{11}\g(s)}\vf_\vp(Z(s))ds\vsp
\qquad+C_{11}\E\dd\int^T_t e^{-C_{11}\g(s)}
\g'(s)\log(|Z(s)|^2+\vp)ds+
e^{-C_{11}\g(t)}-e^{C_{11}\g(T)}.\earr$$
Then, by \eqref{e4.19}, we obtain (because we may assume that $C_{11}>C_7$ and $\vp\le1$)
\begin{equation}\label{e4.28}
\barr{l}
\E[e^{-C_{11}\g(t)}\log(|Z(t)|^2+\vp)]
\le \E[e^{-C_{11}\g(T)}\log(|Z(T)|^2+\vp)]\vsp
\qquad+\vf_\vp(Z(0))
+e^{-C_{11}\g(t)}-e^{C_{11}\g(T)}\vsp\qquad
+C_{11}
\E\dd\int^T_t e^{-C_{11}\g(s)}\g'(s)\log(|Z(s)|^2+1)ds.\earr\end{equation}
On the other hand, we have
\begin{equation}\label{e4.29}
\barr{l}
\E\dd\int^T_0\g'(s)\log(|Z(s)|^2+1)ds\\
\qquad\le C_{13}\E\dd\int^T_0
(\|X_1(t)\|^4
+\|X_2(t)\|^4
+ \|X_1(t)\|^2_{(W^{1,4}(\calo))^3}\vsp
\qquad+\,\|X_2(t)\|^2_{(W^{1,4}(\calo))^3}
+1)(1+\log(|X_1(t)|^2+|X_2(t)|^2+1))dt\vsp
\qquad=:C_{14}<\9.\earr\end{equation}
Here is the argument, to show that indeed $C_{14}<\9$.

We have the interpolation inequality
\begin{equation}\label{e4.28a}
\|u\|_{(W^{1,4}(\calo))^3}
\le C\|u\|^{1-\alpha}_{(H^2(\calo))^3}
\|u\|^\a_{(L^2(\calo))^3},\end{equation}for some $\a\in(0,1).$ (The latter is a consequence of the fact that, for \mbox{$\alpha\in (0,1)$}
suitably chosen, we have $D(A^{\alpha})= (H,D(A))_{\alpha}=(H^{2\alpha}(\calo))^3\cap V \subset
(W^{1,4}(\calo))^{3}$, for all $\alpha\in(1/2,1)$.) This yields
$$\barr{l}
\|X_1(t)\|^2_{(W^{1,4}(\calo))^3}+
\|X_2(t)\|^2_{(W^{1,4}(\calo))^3}\vsp
\qquad
\le C(\|X_1(t)\|^{2(1-\a)}_{(H^2(\calo))^3}
+\|X_2(t)\|^{2(1-\a)}_{(H^2(\calo))^3})
(|X_1(t)|^{2\a}+|X_2(t)|^{2\a}),\vsp\hfill \ff t\in(0,T).\earr$$
Taking into account that for an  arbitrary $\beta\in(0,1)$
$$|\log(|z|^2+1)|\le C_\beta(|z|^\beta+1),\ \ff z,$$
 we get by \eqref{e4.3}, \eqref{e4.4prim}  that $C_{14}<\9$.
Then, by \eqref{e4.28}, we obtain
$$\barr{l}
\E[e^{-C_{11}\g(t)}\log(|Z(t)|^2+\vp)]\vsp
\quad
\le
 \E[e^{-C_{11}\g(T)}\log(|Z(T)|^2+\vp)]
 +e^{-C_{11}\g(t)}
 -e^{-C_{11}\g(T)}
 +\vf_\vp(Z(0))+C_{14}.\earr$$Letting $\vp\to0$, we get \eqref{e4.6a}, as desired.

\begin{remark}\label{r4.1}\rm One might suspect that a controllability result similar to Theo\-rem \ref{t2} remains valid in this case too. However, this requires a forward uniqueness result for the linearized backward stochastic equation corres\-pon\-ding to \eqref{e4.1} which, as in the case of equation \eqref{e3.13}, remains open.
\end{remark}

\begin{remark}\label{r4.2}\rm By inspecting the previous proof, it is clear that Theorem \ref{t3} remains true for any pair of solutions $X_1,X_2$ to the stochastic Navier--Stokes equation
$$\barr{l}
dX-\nu\Delta X dt+(X\cdot\nabla)Xdt=XdW\mbox{ in } (0,T)\times\calo,\vsp
\nabla\cdot X=0,\ X=0\mbox{ on }(0,T)\times\pp\calo.\earr$$
which satisfies condition \eqref{e4.3} (if any). Anyway, Theorem \ref{t3} remains true for linear Oseen--Stokes equations of the form
$$\barr{l}
dX-\nu\D X,dt+((X\cdot\nabla)a+(b\cdot\nabla)X)dt=X\,dW\vsp
\nabla\cdot X=0,\ X=0\ \mbox{ on }(0,T)\times\pp\calo.\earr$$
\end{remark}

\n{\bf Acknowledgements.} Viorel Barbu  was    supported by a grant of the Romanian National Authority for Scientific Research, CNCS-UEFISCDI project PN-II-ID-PCE-2011-3-0027 and Bibos-Research Centre.  Michael R\"ockner's research was supported by DFG through CRC 701.

\end{document}